\documentclass[12pt]{amsart}
\usepackage{times}
\usepackage{amsfonts}
\usepackage{amsthm}
\usepackage{amsmath}
\advance \topmargin by -\headheight
\advance \topmargin by -\headsep
\evensidemargin \oddsidemargin
\newtheorem{theorem}{Theorem}[section]
\newtheorem{lemma}[theorem]{Lemma}
\newtheorem{corollary}[theorem]{Corollary}

\newtheorem{remark}[theorem]{Remark}
\newtheorem{example}[theorem]{Example}

\newtheorem{definition}[theorem]{Definition}

\begin{document}

\title{Strongly $1$-bounded von Neumann algebras}

\author{Kenley Jung}

\address{Department of Mathematics, University of California,
Los Angeles, CA 90095-1555,USA}

\email{kjung@math.ucla.edu}
\thanks{Research supported by an NSF postdoctoral fellowship}

\begin{abstract} Suppose $F$ is a finite set of selfadjoint elements in a 
tracial von Neumann algebra $M$.  For $\alpha >0$, $F$ is $\alpha$-bounded 
if $\mathbb P^{\alpha}(F) < \infty$ where $\mathbb P^{\alpha}$ is the free 
packing $\alpha$-entropy of $F$ introduced in [9].  We say that $M$ is 
strongly $1$-bounded if $M$ has a $1$-bounded finite set of selfadjoint 
generators $F$ such that there exists an $x \in F$ with $\chi(x) > 
-\infty$.  It is shown that if $M$ is strongly $1$-bounded, then any 
finite set of selfadjoint generators $G$ for $M$ is $1$-bounded and 
$\delta_0(G) \leq 1$; consequently, a strongly $1$-bounded von Neumann 
algebra is not isomorphic to an interpolated free group factor and 
$\delta_0$ is an invariant for these algebras.  Examples of strongly 
$1$-bounded von Neumann algebras include (separable) 
$\mathrm{II}_1$-factors which have property 
$\Gamma$, have Cartan subalgebras, are non-prime, or the group von Neumann 
algebras of $SL_n(\mathbb Z)$, $n \geq 3$.  If $M$ and $N$ are strongly 
$1$-bounded and $M \cap N$ is diffuse, then the von Neumann algebra 
generated by $M$ and $N$ is strongly $1$-bounded.  In particular, a free 
product of two strongly $1$-bounded von Neumann algebras with amalgamation 
over a common, diffuse von Neumann subalgebra is strongly $1$-bounded.  
It is also shown that a $\mathrm{II}_1$-factor generated by the normalizer 
of a strongly $1$-bounded von Neumann subalgebra is strongly $1$-bounded.

\end{abstract} \maketitle

\section*{introduction}

Given a finite set of selfadjoint elements $F = \{x_1,\ldots, x_n\}$ in a 
tracial von Neumann algebra $(M, \varphi)$, the $(m,k,\gamma)$-microstate 
space for $F$, $\Gamma(F;m,k,\gamma)$, consists of all $n$-tuples of 
selfadjoint $k\times k$ complex matrices $(a_1, \ldots, a_n)$ such that for $1 \leq p 
\leq m$ and $1 \leq i_1, \ldots, i_p \leq n$

\[ |tr_k(a_{i_1} \cdots a_{i_p}) - \varphi(x_{i_1} \cdots x_{i_p})|< 
\gamma \]

\noindent where $tr_k$ is the normalized tracial state on the $k \times k$ 
matrices. These microstate spaces are subsets of Euclidean spaces and 
hence, Lebesgue volume and packing dimension can be applied to analyze 
them. Voiculescu introduced these notions in [16] and used them to define 
the free entropy $\chi(F)$ of $F$ and the free entropy dimension 
$\delta_0(F)$ of $F$.  $\chi(F)$ is an asymptotic logarithmic volume of 
the microstates of $F$ and $\delta_0(F)$ is an asymptotic 
packing/Minkowski dimension of the microstates of $F$.

\emph{The} issue in microstates theory is the invariance problem for $\delta_0$: if $F$ and $G$ are 
finite sets of selfadjoint elements in $M$, and $F$ and $G$ generate the 
same von Neumann algebra, then is it the case that $\delta_0(F) = 
\delta_0(G)$?  An affirmative answer to this would show the nonisomorphism 
of the free group factors.

[9] studied the microstate spaces with elementary techniques from fractal 
geometry.  This attempt to strengthen the connections between microstate 
theory and geometric measure theory was driven in part by the following two facts, one from free probability, the other from geometric 
measure theory.

On the free probability side all applications in [5] and [17] show 
that von Neumann algebras with certain decomposition properties (Property 
$\Gamma$, Cartan subalgebras, tensor decomposition) satisfy the condition 
that for any finite generating set $F$ of the von Neumann algebra, 
$\delta_0(F)
 \leq 1$.  Assuming the algebra embeds into the ultraproduct of the 
hyperfinite $\mathrm{II}_1$-factor, [7] shows that $1 \leq \delta_0(F)$ so 
that $\delta_0(F)=1$.  Thus, $\delta_0$ is an invariant for such von 
Neumann algebras and their free entropy dimension is $1$.  The free group 
factor $L(F_n)$ has a finite set of selfadjoint generators $X$ such that 
$\delta_0(X) =n$, thus, [5] and [17] show that in fact a free group 
factors cannot have any of these decomposition properties.  
Significantly, these were the first known kind with separable predual 
which are prime or fail to have Cartan subalgebras (Popa shows in [13] 
that the von Neumann algebra on a free group with uncountable many 
generators is prime and has no Cartan subalgebras; but unfortunately, 
these von Neumann algebras are inseparable).

On the geometric measure theory side, Besicovitch classified metric spaces 
with finite Hausdorff $1$-measure (these sets automatically have Hausdorff 
dimension $1$).  His study concluded with the following, fairly complete 
answer: any such space $\Omega$ breaks up into a good and bad part.  The 
good part consists of some Cantor dust and a set which has a tangent at 
almost all of its points; moreover, this latter set can be contained in a 
countable union of rectifiable curves.  The bad part is a totally 
irregular set (all local densities have different upper and lower bounds) 
and no point of the set has a tangent.  For the study of sets with 
nondegenerate Hausdorff $r$-measure with $r>1$ the situation was much more 
complicated and it was some time (about 50 years after Besicovitch's work 
initial work) before some of the basic problems were resolved (see [3] for 
an overview).

The analysis in [5] and [17] decomposes microstate spaces of von Neumann 
algebras with certain properties into microstates of hyperfinite algebras 
and sets of negligible packing entropy. Similarly Besicovitch's work 
decomposes $\Omega$ into rectifiable curves and sets of measure $0$ 
(ignoring the irregular part).  Both approaches express their respective 
problems in terms of well-understood spaces - the injective one in the 
microstate setting, and the real line in the fractal setting.  Given some 
of the already existing connections between the microstates theory and 
geometric measure theory as well as Besicovitch's success in classifying 
sets with finite Hausdorff measure $1$, it seems plausible that the 
invariance problem for $\delta_0$ can be answered affirmatively for finite 
sets with dimension $1$.  More specifically, is it true that if 
$\delta_0(F) =1$, then for any other finite set of selfadjoint generators 
$G$ for $F^{\prime \prime}$, $\delta_0(G) =1$?  Under some additional 
conditions the answer is "yes", and moreover, one can use a decomposition 
argument akin to Besicovitch's where amenability takes the place of 
$[0,1]$.

Suppose $F \subset M$ is a finite set of selfadjoint elements such that 
$\mathbb P^1(F) <\infty$.  This analytic condition on $F$ is called 
$1$-boundedness.  Here $\mathbb P^1(F)$ is a kind of packing $1$-measure 
(this assumption can be likened to the assumption in Besicovitch's 
classification that the set have bounded Hausdorff $1$-measure, though 
strictly speaking, our assumption is stronger).  Assume moreover that $F$ 
contains an element $x$ with finite free entropy.  We show then that any 
other finite set of selfadjoint generators $G$ for $F^{\prime \prime}$ is 
$1$-bounded and that $\delta_0(G)  \leq1$.  The argument uses a microstate 
decomposition relative to a hyperfinite algebra, an idea which appeared in 
qualitative form in [11].  It can be regarded as a Fubini-type theorem 
where, as in Besicovitch's theorem the decomposition breaks up the good 
part of the space into a negligible set of "Cantor dust" (relative 
microstates) and a "rectifiable subset" (hyperfinite microstates).  
Now, when $F^{\prime \prime}$ embeds into the ultraproduct of the 
hyperfinite $\mathrm{II}_1$-factor, then for any such generating $G$ for 
$F^{\prime \prime}$, $\delta_0(G) =1$ (this is a consequence of [7]).  If 
this is not the case, then $\delta_0(G) = -\infty$.  From these facts, it 
follows that $\delta_0$ is indeed an invariant for all von Neumann 
algebras with such a generating set $F$. We say that a von Neumann algebra 
$M$ is strongly $1$-bounded if it has such a generating set $F$.  It is a 
consequence of [5] and [17] that if $M$ has property $\Gamma$, a Cartan 
subalgebra, a nontrivial tensor product decomposition, or if $M$ is a 
group von Neumann algebra of $SL_n(\mathbb Z)$, $n \geq 3$, then $M$ is 
strongly $1$-bounded.

$\delta_0$ is an invariant for strongly $1$-bounded von Neumann algebras 
and their amplifications.  Moreover, if $M$ and $N$ are strongly 
$1$-bounded and $A$ and $B$ are diffuse von Neumann subalgebras of $M$ and 
$N$, respectively, such that $A$ and $B$ generate a strongly $1$-bounded 
von Neumann algebra, then the von Neumann algebra generated by $M$ and $N$ 
is strongly one-bounded and thus, not isomorphic to an interpolated free 
group factor.  This implies, in particular, that if $D = M \cap N$ is 
diffuse, then $M *_D N$ is not isomorphic to a free group factor.

The outline of the paper is as follows.   We start with a brief list of notation followed by the first section which defines  
$\alpha$-bounded sets, remarks on some equivalent formulations, and has a 
short list of examples.  The second section collects some facts on the 
decomposition of microstate spaces relative to a single selfadjoint; it is 
the decomposition of the Fubini/Besicovitch-type described above.  The 
third section states and proves the main result.  The fourth consists of 
nonisomorphism applications to amalgamated free products and other types 
of von Neumann algebras.  The fifth and final section is a generalization 
of [5] and [17].  It will imply, in particular, that if a 
$\mathrm{II}_1$-factor can be generated by the normalizer of a 
strongly $1$-bounded von Neumann subalgebra, then $N$ is strongly $1$-bounded. 
It follows that an interpolated free group factor $L(F_r)$ cannot 
be isomorphic to the crossed product of a strongly $1$-bounded von Neumann 
algebra with a group action.  

\section*{notation}

Throughout $M$ denotes a tracial von Neumann algebra with separable 
predual.  For any $k, n \in \mathbb N, R>0$, $M^{sa}_k(\mathbb C)$ 
denotes the $k \times k$ complex selfadjoint matrices, $(M^{sa}_k(\mathbb 
C))^n$ denotes $n$-tuples of entries in $M^{sa}_k(\mathbb C)$ and 
$(M^{sa}_k(\mathbb C))_R$ denotes the elements of $M^{sa}_k(\mathbb C)$ 
with operator norms no greater than $R$.  For $\xi = (\xi_1,\ldots, \xi_n) 
\in (M^{sa}_k(\mathbb C))^n$ we write $|\xi|_2 = (\Sigma_{i=1}^n 
(tr_k(\xi_i^2)))^{\frac{1}{2}}$ where $tr_k$ is the tracial state on 
$M_k(\mathbb C)$, the $k \times k$ matrices, and for a $k \times k$ 
unitary $u$, $u\xi u^* = (u\xi_1 u^*, \ldots, u\xi_n u^*)$.

\section{$\alpha$-bounded sets}

In this section assume $F$ is a finite set of selfadjoint elements of 
$M$.  Recall the definitions of $\mathbb K_{\epsilon, \infty}(F)$ and 
$\mathbb P_{\epsilon, \infty}(F)$ introduced in [11].  We have the 
following definition:

\begin{definition} For $\alpha >0$, $F$ is said to be $\alpha$-bounded if 
for some $K, \epsilon_0 >0$ and any $\epsilon_0 > \epsilon >0,$

\[ \mathbb K_{\epsilon, \infty}(F) \leq  \alpha  \cdot |\log \epsilon|+ K.
\]
\end{definition} 

\begin{remark}It is immediate from [8] and Lemma 2.2 of [11], that if $F$ 
is $\alpha$-bounded, then $\delta_0(F) \leq \alpha$.  \end{remark}

Recall that cutoff constants for the operator norm were used in the 
definition of $\mathbb K_{\epsilon}(F)$ and $\mathbb P_{\epsilon}(F)$ and 
that $\epsilon$ quantities $\mathbb K_{\epsilon, R}(F)$ and $\mathbb 
P_{\epsilon, R}(F)$ were introduced in [8] where the microstates spaces 
have cutoff constants.

\begin{lemma} For any $\epsilon >0$ and $R \geq \max_{x \in F} \{\|x\|\}$ 
we have

\[\mathbb P_{4\epsilon, \infty}(F) \leq \mathbb K_{2\epsilon, \infty}(F) 
\leq \mathbb K_{\epsilon, R}(F) \leq \mathbb K_{\epsilon}(F) \leq \mathbb 
K_{\epsilon, \infty}(F) \leq \mathbb P_{\frac{\epsilon}{2},\infty}(F). \] 
\end{lemma}

\begin{proof} By Lemma 2.1 of [11] $\mathbb K_{ 2\epsilon, \infty}(F) \leq 
\mathbb K_{\epsilon, R}(F) \leq \mathbb K_{\epsilon}(F)$.  The rest of the statement
follows from the fact that for any metric space $\Omega$ and $\epsilon 
>0$, $P_{\epsilon}(\Omega) \geq K_{2 \epsilon}(\Omega) \geq P_{4 
\epsilon}(\Omega)$ where $P_{\epsilon}(\Omega)$ is the maximum number in a 
collection of disjoint open $\epsilon$-ball with centers in $\Omega$ and 
$K_{\epsilon}(\Omega)$ is the minimum number of open $\epsilon$-balls 
required to cover $\Omega$. \end{proof}

The free packing $\alpha$-entropy of $F$ was defined in [9] as $\mathbb 
P^{\alpha}(F) = \sup_{R >0} \mathbb P^{\alpha}_R(F)$ where $\mathbb 
P^{\alpha}_R(F) = \limsup_{\epsilon \rightarrow 0} \mathbb P_{\epsilon, 
R}(F) + \alpha \log 2 \epsilon$.  The free packing $\alpha$-entropy has 
the same relationship to $\delta_0$ that Hausdorff measure (free Hausdorff 
entropy) has to Hausdorff dimension (free Hausdorff dimension).  From 
Lemma 1.3 and Lemma 3.11 of [9] we have:

\begin{corollary} Suppose $F = \{x_1,\ldots, x_n\}$, $\{s_1, \ldots, 
s_n\}$ is a semicircular family free with respect to $F$, and $R \geq 
\max_{x \in F} \{\|x\|\}$.  The following conditions are equivalent: 
\begin{itemize} \item $F$ is $\alpha$-bounded. \item There exist $K, 
\epsilon_0 >0$ such that for all $\epsilon_0 > \epsilon >0$, $\mathbb 
K_{\epsilon, R}(F) \leq \alpha \cdot |\log \epsilon|+K$. \item There exist 
$K, \epsilon_0 >0$ such that for all $\epsilon_0 > \epsilon >0$, $\mathbb 
P_{\epsilon, \infty}(F) \leq \alpha \cdot \log \epsilon| + K$. \item There 
exist $K, \epsilon_0 >0$ such that for all $\epsilon_0 > \epsilon >0$, 
$\mathbb P_{\epsilon, R}(F) \leq \alpha \cdot |\log \epsilon|+K$. \item 
$\mathbb P^{\alpha}(F) < \infty$. \item $\limsup_{\epsilon \rightarrow 0} 
\left (\chi(x_1 + \epsilon s_1, \ldots, x_n + \epsilon s_n: s_1,\ldots, 
s_n) + (n - \alpha) |\log \epsilon| \right ) < \infty$. \end{itemize} 
\end{corollary}

\begin{example} Suppose $F^{\prime \prime}$ can be generated by a sequence 
of Haar unitaries $\langle u_j \rangle_{j=1}^{\infty}$ such that for each 
$j$, $u_{j+1}^* u_j u_{j+1} \in \{u_1,\ldots, u_j\}^{\prime \prime}$.  By 
[5] $F$ is $1$-bounded.  Indeed, it is a consequence of Lemma 5.1 of [5] 
that when $F$ consists of contractions and $\epsilon \in (0,1)$, $\mathbb 
K_{\epsilon}(F) \leq |\log \epsilon| + 100 + \frac{1}{2} \log (\#F)$, and 
this condition clearly implies that $F$ is $1$-bounded for general $F$.

This class of von Neumann algebras considered by Ge and Shen include the 
following cases: 1) $F$ generates a von Neumann algebra with a regular, 
diffuse, hyperfinite von Neumann subalgebra; 2) $F$ generates a von 
Neumann algebra of the form $A \otimes B$ where $A$ and $B$ are diffuse; 
3) $F$ generates the group von Neumann algebras obtained from $SL_n(\mathbb Z)$ for $n \geq 3$.  In Section 5 we will generalize this example.
\end{example}

\begin{example} Suppose $F^{\prime \prime} =M$ and $M^{\omega}$ is the 
ultraproduct of $M$ for some nontrivial $\omega$.  If $A \subset 
M^{\omega}$ is a diffuse hyperfinite subalgebra such that the von Neumann 
algebra generated by the normalizer of $A$ contains $B$, then using the 
free entropy formulation in Corollary 1.4, it follows from Theorem 7.3 of 
[17] that $F$ is $1$-bounded.  In particular, if $M$ has property 
$\Gamma$, then $F$ is $1$-bounded. \end{example}

\begin{example} Suppose $F$ is a finite generating set of the type 
considered in [10] for the interpolated free group factor $L(F_r)$, $r >1$, 
of Dykema and Radulescu ([2] and [14]).  It follow that $F$ is 
$r$-bounded.  By Theorem 7.4 of [9] this also holds if $F$ can be broken 
up into a free collection of subsets, each of which generates a finite 
dimensional algebra. \end{example}

\section{Microstates Relative to a Single Selfadjoint}

In this section we want to find packing entropy estimates with respect to 
finite sets of selfadjoints of the form $\{x\} \cup F$ where $\chi(x) > 
-\infty$.  These will come from the microstates of $F$ relative to a fixed 
sequence of microstates for $x$.  Much of this will be a regurgitation of 
the material in [11] but we do so for completeness and because specific 
properties of such sets will be exploited to give more quantitative 
estimate.  For the remainder of this section $F \subset M$ is a finite set 
of selfadjoint elements and $x \in M$ is selfadjoint.  Fix $R >0$ such 
that $R$ is greater than the operator norm of any element in $\{x\} \cup 
F$.  It is easy to see that there exists a sequence $\langle x_k 
\rangle_{k=1}^{\infty}$ such that for each $k$ $x_k \in M^{sa}_k(\mathbb 
C)$, $\|x_k \| <R$, and for any $m \in \mathbb N$ and $\gamma >0$ $x_k \in 
\Gamma(x;m,k,\gamma)$ for sufficiently large $k$.  Fix this sequence 
$\langle x_k \rangle_{k=1}^{\infty}$.  Recall from [11] the microstate 
spaces $\Xi(F;m,k,\gamma)$ for $F$ relative to $\langle x_k 
\rangle_{k=1}^{\infty}$.

\[ \Xi(F;m,k,\gamma) = \{\xi : (x_k, \xi) \in \Gamma(\{x\} \cup F;m,
k,\gamma) \}.\]

\noindent Define successively for $\epsilon >0,$

\[ \mathbb K_{\epsilon}(\Xi(F;m,\gamma)) = \limsup_{k \rightarrow
\infty} k^{-2} \cdot \log K_{\epsilon}(\Xi(F;m,k,\gamma)),\]

\[ \mathbb K_{\epsilon}(\Xi(F)) = \inf \{ \mathbb K_{\epsilon}(F;m,
\gamma) : m \in \mathbb N, \gamma >0\}, \]

\noindent where the packing quantities are taken with respect to
$|\cdot |_2.$ In a similar fashion, we define $\mathbb
P_{\epsilon}(\Xi(F))$ by replacing the $K_{\epsilon}$ above with 
$P_{\epsilon}.$  We will also use the notation $\Xi_R(F;m,k,\gamma)$ and 
$K_{\epsilon,R}(\Xi(F;m,\gamma))$, $P_{\epsilon,R}(\Xi(F;m,\gamma))$ 
for the quantities and sets where the cutoff constant $R$ is used (so 
these are the relative microstates and associated quantities where the 
operator norms of the entries are all less than $R$).

For a finite set of selfadjoint elements $X$ write 
$\underline{\chi}(X)$ for the quantity obtained by replacing the 
$\limsup_{k\rightarrow \infty}$ in the definition of $\chi(X)$ with a 
$\liminf_{k\rightarrow \infty}$.  $\chi(X) \geq \underline{\chi}(X)$ (when 
equality occurs $X$ is said to be regular, see [18]).  We also write 
$\underline{\mathbb H}^{\alpha}(X)$ and $\underline{\mathbb P}_t(X)$ for 
the quantities obtained by replacing the $\limsup_{k\rightarrow \infty}$ 
in the definitions of $\mathbb H^{\alpha}(X)$ and $\mathbb P_t(X)$ with a 
$\liminf_{k \rightarrow \infty}$.

\begin{lemma} If $\{x\} \cup F$ is an $\alpha$-bounded set and $\chi(x) > 
-\infty$, then there exist constants $C, \epsilon_1 >0$ dependent on 
$\{x\} \cup F$ such that for $\epsilon_1 > t >0$

\[ (\alpha - 1) |\log t| +C \geq \mathbb K_t(\Xi(F)).
\]
\end{lemma}

\begin{proof} By Proposition of 4.5 of [16] and Lemma 3.7 of [9] we have 
$\underline{\mathbb H}^1(x) = \underline{\chi}(x) + \frac{1}{2} \log 
(\frac{2}{\pi e}) =\chi(x) + \frac{1}{2} \log (\frac{2}{\pi e}) > - 
\infty$.  Thus, there exists an $r >0$ such that for all $r > t >0$, 
$\underline{\mathbb H}^1_{t}(x) \geq c$ where $c = \chi(x) + \frac{1}{2} 
\log(\frac{2}{\pi e}) -1.$ It is easy to see that for such $t$,

\[ c \leq \underline{\mathbb H}^1_t(x) \leq \underline{\mathbb P}_t(x) + 
\log 4 t.
\]

\noindent Thus, for $r > t >0,$ $c - \log 4 + |\log t| < 
\underline{\mathbb P}_t(x)$.  Also, $\{x\} \cup F$ is $\alpha$-bounded; 
let $\epsilon_0 >0$ and $K$ be as in the definition of 
$\alpha$-boundedness.

Suppose $0 < t < \min\{r, \epsilon_0\}$.  There exist $m_1 \in \mathbb N$ 
and $\gamma_1 >0$ such that for all $m > m_1$ and $0 < \gamma < \gamma_1$,

\begin{eqnarray} \mathbb P_{t}(F;m,\gamma) \leq r \cdot |\log 
t| +K.
\end{eqnarray}

\noindent Clearly for all $m \in \mathbb N$ and $\gamma >0$,

\begin{eqnarray} c -\log 4 + |\log t| \leq \underline{\mathbb 
P}_t(x;m,\gamma).
\end{eqnarray}

\noindent By [7] and [11], there exist $m_2 \in \mathbb N$ and $\gamma_2 
>0$ such that if $a, b \in \Gamma(x,m_2,k,\gamma_2)$, then there exists a 
$k\times k$ unitary $u$ such that $|uau^* -b|_2 < \frac{t}{100}$. 
Combining this with (2) it follows that there exist for $m \geq m_1+ m_2$, 
$0 < \gamma < \min\{\gamma_1, \gamma_2\}$, and $k$ sufficiently large, 
unitaries $\langle v_{\lambda k} \rangle_{ \lambda \in \lambda_k}$ such 
that the balls of $\Gamma(x;m,k,\gamma)$ with centers $\langle v_{\lambda 
k} x_k v^*_{\lambda k} \rangle_{\lambda \in \Lambda_k}$ and radii 
$\frac{99t}{100}$ form a collection of disjoint subsets of 
$\Gamma(x;m,k,\gamma)$ and $\liminf_{k \rightarrow \infty} k^{-2} \cdot 
\log \# \Lambda_k = \underline{\mathbb P}_{t}(x;m,\gamma)$. For each $m 
\geq m_1 + m_2$, $\min\{\gamma_1, \gamma_2\} > \gamma >0$, and 
sufficiently large $k$ pick a subset $\langle \xi_{jk} \rangle_{j \in 
J_k}$ of $\Xi(F;m,k,\gamma)$ of maximal cardinality with respect to the 
condition that the $\epsilon$-balls of $\Xi(F;m,k,\gamma)$ with centers 
$\xi_{jk}$ are disjoint.  It is easy to see that the balls of 
$\Gamma(\{x\} \cup F;m,k,\gamma)$ with centers

\[ \langle (v_{\lambda k} x_k v_{\lambda k}^*, v_{\lambda k} \xi_{jk} 
v_{\lambda k}^*) \rangle_{(\lambda, j) \in \Lambda_k \times J_k}
\]

\noindent and radii $\frac{99t}{100}$ are a pairwise dijoint.  So, using 
(1) and (2) we have for $m \geq m_1 +m_2$ and $0 < \gamma < 
\min\{\gamma_1,\gamma_2\}$

\begin{eqnarray*} \alpha |\log t| + K & \geq & \mathbb P_{t2^{-1}} 
(\{x\} \cup F;m, \gamma) \\ & \geq & \limsup_{k\rightarrow \infty} k^{-2} 
\cdot \log (\# \Lambda_k \cdot \#J_k)\\ &\geq & \liminf_{k \rightarrow 
\infty} k^{-2} \cdot \log \#\Lambda_k + \limsup_{k \rightarrow \infty} 
k^{-2} \cdot \log P_t (\Xi(F;m,k, \gamma)) \\ & \geq & \underline{\mathbb 
P}_t(x;m,\gamma)+ \limsup_{k \rightarrow \infty} k^{-2} \cdot \log 
P_t(\Xi(F; m,k, \gamma)) \\ & \geq & c -\log 4 + |\log t|
+ \mathbb P_t(\Xi(F)).
\end{eqnarray*}

\noindent Grouping the constants together on one side we have for all 
$\min\{r, \epsilon_0\} > t > 0$, 

\[ (\alpha-1) |\log t| + 
K - c + \log 4 \geq \mathbb P_{t}(\Xi(F)) \geq \mathbb K_{2t}(\Xi(F)).\]  

\end{proof}

\begin{lemma} Suppose $B = (\Sigma_{y \in \{x\} \cup F} 
\|x\|_2^2)^{\frac{1}{2}}$.  There exists $L, \epsilon_2 >0$ independent of $\{x\} \cup F$ such that for $0 < \epsilon < \epsilon_2$,

\[ \mathbb K_{\epsilon, R}(\{x\} \cup F) \leq \log((4B+6)L) 
+ |\log \epsilon| + 
\mathbb K_{\frac{\epsilon}{4B+ 6}, R}(\Xi(F)).
\]
\end{lemma}

\begin{proof} By [15] or [19] there are $L, \epsilon_0 >0$ such that for 
$\epsilon_0 > \epsilon >0$ and for any $k \in \mathbb N$ there exists an 
$\epsilon$-net for $U_k$ with respect to the quotient metric induced by $| 
\cdot|_{\infty}$ with cardinality no greater than 
$\left(\frac{L}{\epsilon}\right)^{k^2}.$ Suppose $m \in \mathbb N$ and 
$\gamma >0$.  Observe that there exists $\epsilon > r >0$ so small that 
for any $k$, if $(\xi, \eta) \in \Gamma_R(\{x\} \cup F;m,k,\gamma/2)$, 
$|(\xi, \eta) - (a, b)|_2 < r$, and all the entries of $(a,b)$ have 
operator norms less than or equal to $R$, then $(x, a) \in \Gamma_R(\{x\} 
\cup F;m,k,\gamma)$. There also exist $m_1 \in \mathbb N$ and $\gamma_1 
>0$ such that if $y,z \in \Gamma(x;m_1,k,\gamma_1),$ then there exists a 
$k \times k$ unitary $u$ satisfying $|uy u^* -z|_2 < r.$ Finally, we can 
find $m_2$ and $\gamma_2$ such that for any $(\xi, \eta) \in \Gamma(\{x\} 
\cup F;m_2, k,\gamma_2)$, $|(\xi, \eta)|_2 < B +1$. Set $m_3 =m
+ m_1 + m_2$ and 
$\gamma_3 = \min \{ \gamma/2, \gamma_1, \gamma_2).$

For each $k$ find an $\epsilon$-net $\langle \eta_{jk} \rangle_{j \in 
J_k}$ for $\Xi_R(F;m,k,\gamma)$ with respect to $| \cdot |_2$ of minimum 
cardinality.  Find for each $k$ a set of unitaries $\langle u_{gk} 
\rangle_{g \in G_k}$ such that they form  
$\epsilon$-net with respect to the operator norm and such that

\[ \# G_k \leq \left (\frac{L}{\epsilon} \right)^{k^2}.\]

\noindent Consider

\[ \langle (u_{gk} x_k u_{gk}^*, u_{gk}\eta_{jk}u_{gk}^*)
\rangle_{(g,j) \in G_k \times J_k}.\]

\noindent I claim that this set is a $(4B+6) \epsilon$-net for $\Gamma_R(\{x\} 
\cup F;m_2,k,\gamma_2).$

To see this suppose $(\xi, \eta) \in \Gamma_R(\{x\} \cup F;m_3, k, 
\gamma_3).$ By the selection of $m_1$ and $\gamma_1$ there exists a $u \in 
U_k$ such that $|u^* x_k u - \xi|_2 < r_.$ Taking into account the 
stipulation on $r$ this implies that ($u^* x_k u, \eta) \in \Gamma_R(\{x\} 
\cup F;m,k,\gamma) \Longleftrightarrow (x_k, u \eta u^*) \in 
\Gamma_R(\{x\} \cup F;m,k,\gamma),$ whence $u \eta u^* \in 
\Xi_R(F;m,k,\gamma).$ There exists an $g \in G_k$ such that $\|u - u_{gk} 
\| 
< \epsilon.$ Consequently,

\[ |u^*_{gk}x_k u_{gk} - \xi|_2 \leq 2(B+1)\epsilon  + |u^* x_k u - \xi|_2 
\leq (2 B+3) \epsilon. \]

\noindent $u \eta u^* \in \Xi_R(F;m,k,\gamma)$ so there exists a $j \in
J_k$ such that $|\eta_{jk} - u \eta u^* |_2 < \epsilon.$ Now,

\[ |u_{gk}^* \eta_{j k} u_{gk} - \eta|_2 < 2 (B +1)  \epsilon + |u^* \eta_{jk} u 
- \eta|_2 < (2B+3) \epsilon. \]

\noindent $|(u_{gk}^* x_k u_{gk}, u_{gk}^* \eta_{jk} u_{gk}) - (\xi, \eta)|_2 < (4B+6) 
\epsilon$ as desired.

It follows that

\begin{eqnarray*} \mathbb K_{(4 B +6)\epsilon, R}(\{x\} 
\cup F;m_3, \gamma_3) & \leq & \limsup_{k 
\rightarrow \infty} k^{-2} \cdot \log( \# G_k 
\cdot \# J_k ) \\ & \leq & \log L + |\log 
\epsilon| + \limsup_{k \rightarrow \infty} 
k^{-2} \cdot \log 
K_{\epsilon}(\Xi_R(F;m,k,\gamma)).  
\end{eqnarray*}

\noindent Given $0 < \epsilon < \epsilon_0$ and 
any $m \in \mathbb N$ and $\gamma >0$ we 
produced $m_3 \in\mathbb N$ and $\gamma_3 >0$ so 
that the above inequality holds. Thus

\[ \mathbb K_{(4 B+6)\epsilon,R} (\{x\} \cup F) \leq \log 
L +|\log\epsilon| + \mathbb 
K_{\epsilon, R}(\Xi(F)). \]

\noindent This clearly implies that for $0 < \epsilon < (4B+6) \epsilon_0$

\[ \mathbb K_{\epsilon, R}(\{x\} \cup F) \leq \log ((4B+6)L) 
+ 
|\log \epsilon| + \mathbb 
K_{\frac{\epsilon}{4B+6}, R}(\Xi(F)).
\]
\end{proof}

\section{Strongly $1$-bounded von Neumann algebras}

We now come to our main result which says that if $\{x\} \cup
F$ is a $1$-bounded set of selfadjoint generators for $M$ such that
$\chi(x) > -\infty$, then any other finite set of selfadjoint generators
$G$ for $M$ is $1$-bounded.  First a lemma:

\begin{lemma} If $X$ and $Y$ are finite sets of selfadjoint elements such 
that $Y \subset X^{\prime \prime}$, then for any $R, \epsilon >0$, $\mathbb K_{\epsilon, R}(X) \leq \mathbb K_{\epsilon, R}(X \cup Y)$.  In particular, if $X \cup Y$ is $r$-bounded, then $X$ is $r$-bounded. \end{lemma}

\begin{proof} It suffices to prove the first statement.  This is a 
repetition of Lemma 3.6 in [9].  Suppose $R$ exceeds the operator norms of 
any of the elements in $X \cup Y$.  Given $m \in \mathbb N$ and $\epsilon, 
\gamma >0$ there exist $m_1 \in \mathbb N,$ $R, \gamma_1>0$ and a 
$\#Y$-tuple $f$ of polynomials in $n$ noncommutative variables such that 
if $\xi \in \Gamma_R(X;m_1,k,\gamma_1)$ then

\[ (\xi, f(\xi)) \in
\Gamma_R(X \cup Y ;m,k,\gamma).\]

\noindent For each $k$ this map from
$\Gamma_1(X;m_1,k,\gamma_1)$ to
$\Gamma_R(X \cup Y;m,k,\gamma)$ defined by sending
$\xi$ to $(\xi, f(\xi))$ increases distances with respect to $|\cdot|_2.$
Hence

\[ \mathbb K_{\epsilon, R}(X;m_1, \gamma_1) \leq \mathbb K_{\epsilon, R}(X 
\cup Y;m,\gamma).\]

\noindent This being true for any $m,\gamma, \epsilon$ the result 
follows from Corollary 1.4.  
\end{proof} 

\begin{theorem} If $\{x\} \cup F$ is a $1$-bounded finite set of 
selfadjoint generators for $M$ such that $\chi(x) > -\infty$, then for any 
other finite set of selfadjoint generators $G$ for $M$, $G$ is a 
$1$-bounded set.  In particular, for such $G$, $\delta_0(G) \leq 1$. 
\end{theorem}

\begin{proof} Suppose $G$ is a finite set of selfadjoint generators for 
$M$.  By Lemma 3.1 in order to show that $G$ is $1$-bounded, it suffices 
to show that $\{x\} \cup F \cup G$ is $1$-bounded.

Fix $R$ such that $R$ exceeds the operator norms of any of the elements in 
$\{x\} \cup F \cup G$.  Find a sequence $\langle x_k 
\rangle_{k=1}^{\infty}$ with the associated $C, \epsilon_1 >0$ as 
constructed in Lemma 2.1 satisfying the covering bound for the relative 
microstates $\Xi(F;m,k,\gamma)$.  Applying Lemma 2.2 to the set $F \cup 
G$, if $B = (\Sigma_{y \in \{x\} \cup F \cup G} \|y\|_2^2)^{\frac{1}{2}}$, 
then there exist constants $K, \epsilon_2 >0$ with $K$ dependent only on 
$B$ such that for all $\epsilon_2 > \epsilon >0$,

\begin{eqnarray} \mathbb K_{\epsilon, R}(\{x\} \cup F \cup G) \leq K + 
|\log \epsilon| + \mathbb K_{\frac{\epsilon}{4B+6},R}(\Xi(F \cup G)).
\end{eqnarray}

\noindent Set $D = 4B+6$ and suppose $\epsilon_2 >\epsilon >0$.  There exists an 
$\#G$-tuple of polynomials in $\#F+1$ noncommuting variables, $\Phi$, such 
that $|\Phi(x,F) - G|_2 < \epsilon(10D)^{-1}$. Moreover, 
  regarding $\Phi$ as a map from $((M^{sa}_k(\mathbb C))_R)^{\#F+1} 
  \rightarrow (M^{sa}_k(\mathbb C))^{\#G}$, $\Phi$ has a bounded Lipschitz 
  constant $L_k$ with respect to the $| \cdot |_2$-norm and $\sup_{k \in \mathbb N} L_k < \infty$.  Choose $L$ greater 
  than $\sup_{k \in \mathbb N} L_k$ and greater than $\max \{\epsilon^{-1}, 1\}$.

The selection of $C, \epsilon_1 >0$ implies that for $2^{-1} \cdot 
\min\{\epsilon (10(DL))^{-1}, \epsilon_1\} > t >0$,

\begin{eqnarray*} C \geq & \mathbb K_t (\Xi(F)). \end{eqnarray*}

\noindent So, there exist $m_1 \in \mathbb N$ and $\gamma_1 >0$ such that 
for $m \geq m_1$ and $\gamma_1 > \gamma >0$,

\begin{eqnarray} C +1 \geq \limsup_{k 
\rightarrow \infty} k^{-2} \cdot \log K_t (\Xi_R(F;m,k,\gamma)). 
\end{eqnarray}

Now there clearly exist $m_2 \in \mathbb N$ and $\gamma_2 >0$ such that if 
$m > m_2$ and $\gamma_2 > \gamma >0$, then for any $k$, $(\xi, \eta) \in 
\Xi( F \cup G;m,k,\gamma) \Rightarrow |\Phi(x_k,\xi) - \eta|_2 \leq 
\epsilon (10D)^{-1}$.  Suppose $m \in \mathbb N$ and $\gamma >0$ with $m > 
m_1 +m_2$ and $\min\{\gamma_1, \gamma_2\} > \gamma$.  By (4) we can find 
for $k$ sufficiently large, an $t$-net $\langle \xi_{jk} \rangle_{j \in 
J_k}$ for $\Xi_R(F;m,k,\gamma)$ which satisfies

\begin{eqnarray} \#J_k \leq e^{(C+1)k^2}. \end{eqnarray}

\noindent For each such $k$ sufficiently large consider the set $\langle 
(\xi_{jk}, \Phi(x_k,\xi_{jk})) \rangle_{j \in J_k}$.  I claim that this 
set is an $\epsilon D^{-1}$-cover for $\Xi_R(F \cup G;m,k,\gamma)$.  
Indeed, suppose $(\xi, \eta) \in \Xi_R(F \cup G;m,k,\gamma)$.  Then by 
definition, $\xi \in \Xi_R(F;m,k,\gamma)$ so that there exists some $j_0 
\in J_k$ with $|\xi
- \xi_{j_0}|_2 < \epsilon (10 D L)^{-1}$.  Since both $\xi$ and 
  $\xi_{j_0k}$
  are $\#F$-tuples of operators with norms no greater than $R$, 
$|\Phi(x_k, \xi) -\Phi(x_k, \xi_{j_0k})|_2 < L 
  \cdot t \leq \epsilon (10D)^{-1}.$ On the 
other hand, 
  $|\Phi(x_k,\xi) - \eta|_2 \leq \epsilon (10 D)^{-1}$ so that 
$|\Phi(x_k, \xi_{j_0 k}) - \eta|_2
  < \epsilon (5D)^{-1}.$ Thus, $|(\xi, \eta) - 
(\xi_{j_0k},\Phi(x_k,\xi_{j_0k}))|_2 
<  t + \epsilon (5D)^{-1} < \epsilon D^{-1}.$

By the preceding paragraph and (5) we conclude that for $m > m_1 + m_2$, 
$\min\{\gamma_1, \gamma_2\} > 
\gamma >0$ and $k$ sufficiently large,

\[k^{-2} \cdot \log \left[K_{\frac{\epsilon}{4B+ 6}, R}(\Xi(F \cup 
G;m,k,\gamma)) 
\right] \leq C+1. \]

\noindent Taking a $\limsup_{k \rightarrow \infty}$ on both sides yields

\begin{eqnarray} \mathbb K_{\frac{\epsilon}{4B+ 6}, R}(\Xi(F \cup G))  
\leq  
\mathbb K_{\frac{\epsilon}{4B+ 6}, R}(\Xi(F \cup G;m,\gamma)) \leq  C + 1 
\end{eqnarray}

\noindent Stuffing (6) into (3) yields for all $0 < \epsilon < \epsilon_2$

\[ \mathbb K_{\epsilon, R}(\{x\} \cup F \cup G) \leq K + C + 1 + |\log 
\epsilon|.
\]

\noindent By Corollary 1.4, $\{x\} \cup F \cup G$ is $1$-bounded, and 
thus, by Lemma 3.1 so is $G$.
\end{proof}

$1$-boundedness is a condition on a finite set and its microstate spaces, 
and makes no direct reference to the generated algebra.  The point of 
Theorem 3.2 is that $1$-boundedness of an appropriate finite set imposes 
$1$-boundedness on any other generating set of the von Neumann algebra of 
the initial set.  In this way, $1$-boundedness is a property of the set 
which propagates to a property of the generated von Neumann algebra.

In view of Theorem 3.2, we make the following definition:

\begin{definition} $M$ is strongly $1$-bounded if $M$ has a $1$-bounded 
finite set of selfadjoint generators $F$ such that $F$ contains an 
element $x$ with $\chi(x) > -\infty$. \end{definition}

\begin{remark} Strong $1$-boundedness of a tracial von Neumann algebra 
requires firstly, that the algebra have a finite set of selfadjoint 
generators with one element having finite free entropy and secondly, that 
the free packing $1$-entropy of such a generating set be finite.  
The first condition is equivalent to having the von Neumann algebra 
possess a finite generating set.  It is unknown whether an arbitrary von 
Neumann algebra (with separable predual) can be generated by a finite set 
of selfadjoint elements. However, it is well-known that the examples 
considered in Section 1 are finitely generated when factoriality is 
imposed.  Suppose $M$ is a $\mathrm{II}_1$-factor (with separable predual, 
as is always tacitly assumed).  By [4] if $M$ has property $\Gamma$, then 
$M$ can be generated by a finite number of elements; thus, by Example 1.6 
it follows that $M$ is strongly $1$-bounded.  Similarly, by [6], if $M$ 
has a regular, diffuse, hyperfinite von Neumann subalgebra or is 
non-prime, then $M$ has a finite set of selfadjoint generators and by 
Example 1.5, it follows that $M$ is strongly $1$-bounded. \end{remark}

By Theorem 3.2 we have:

\begin{corollary} If $M$ is strongly $1$-bounded, then for any finite set 
of selfadjoint generators $X$ of $M$, $\delta_0(X) \leq 1$.  Moreover, if 
$M$ is embeddable into the ultraproduct of the hyperfinite 
$\mathrm{II}_1$-factor, then for any finite set of selfadjoint generators 
$X$ of $M$, $\delta_0(X) =1$.  Thus, $\delta_0$ is an invariant for these 
algebras.  \end{corollary}

We have by Lemma 5.2 of [12]

\begin{corollary} If $M$ is strongly $1$-bounded and $\alpha >0$, then for 
any finite set of selfadjoint generators $X$ for $M_{\alpha}$, the 
amplification of $M$ by $\alpha$, $\delta_0(X) \leq 1$. \end{corollary}

Recall that Dykema and Radulescu ([2], [14]) defined a family of von 
Neumann algebras, $L(F_r)$, $1 < r \leq \infty$ such that for integer 
values $r \in \mathbb N \cup \{\infty\}$, $L(F_r)$ coincides with the free 
group factor on $r$ generators.  These von Neumann algebras are called the 
interpolated free group factors.  Voiculescu was the first to show ([16]) 
that there exists a finite set of generators $X$ for $L(F_r)$ such that 
$\delta_0(X) =r$.  Thus, by Corollary 3.5 if $M$ is strongly $1$-bounded, 
then $M$ cannot be isomorphic to $L(F_r)$ for $1 < r < \infty$.  We 
include $\epsilon$ more by shows that $M$ cannot be isomorphic to 
$L(F_{\infty})$:

\begin{lemma} If $1< s \in \mathbb N \cup \{\infty\}$ and $\langle 
M_i \rangle_{i=1}^{s}$ is a sequence of finitely generated diffuse von Neumann subalgebras of 
$M$ such that each $M_i$ embeds into the ultraproduct of the hyperfinite 
$\mathrm{II}_1$-factor, then $*_{i=1}^s M_i$ cannot be strongly 
$1$-bounded.\end{lemma}

\begin{proof} Suppose to the contrary that $*_{i=1}^s M_i$ can be strongly 
$1$-bounded.  By definition there exists a finite set of selfadjoint 
generators $X$ for $*_{i=1}^s M_i$.  There exist finite sets of 
selfadjoint elements, $F_j$ such that $F_j^{\prime \prime} = M_j \subset 
*_{i=1}^s M_i$, $j =1, 2$.  By [7] and [17] it follows that $2 \leq 
\delta_0(F_1 \cup F_2)$.  The embeddability assumption 
on the $M_i$ and the asymptotic freeness results of [18] imply 
that $\delta_0(F_1 \cup F_2) \leq \delta_0(X \cup F_1 \cup F_2)$.  On the 
other hand, $X \cup F_1 \cup F_2$ is a finite set of generators for 
$*_{i=1}^s M_i$ so Theorem 3.2 implies that $X \cup F_1 \cup F_2$ is 
$1$-bounded; consequently by Remark 1.2, $\delta_0(X \cup F_1 \cup F_2) 
\leq 1$.  Putting this together we have

\[ 2 \leq \delta_0(F_1 \cup F_2) \leq \delta_0(X \cup F_1 \cup F_2) \leq 
1. \]

\noindent This is preposterous.  $*_{i=1}^s M_i$ cannot be strongly $1$-bounded. 
\end{proof}

In [1], Nate Brown constructed finite sets $X$ of selfadjoint elements 
such that $\chi(X) > -\infty$, $\#X >1$, and $X^{\prime \prime} \neq 
L(F_r)$.  In particular by [17] $\delta_0(X) = \#X \neq 1$ for such $X$.  
Combining all this with Corollary 3.5 and Lemma 3.7 above we have

\begin{corollary} If $M$ is strongly $1$-bounded, then 
$M$ is neither isomorphic to an interpolated free group 
factor nor to the free perturbation algebras considered in [1].  
\end{corollary}

\section{An application to amalgamated free products}

In this section we want to produce some other examples of 
$1$-bounded sets.  In particular we 
show (see Corollary 4.2 and Corollary 4.4) that if $A$ 
and $B$ are von Neumann algebras which are nonprime, have 
Cartan subalgebras, or have property $\Gamma$ and $D$ is a 
diffuse subalgebra of $A$ and $B$, then $A *_D B$ is not 
isomorphic to an interpolated free group factor.   
The proof will again rest on the relative hyperfinite 
decomposition results in Section 2.

If $M$ and $N$ are von Neumann algebras acting on the same Hilbert space, 
then we denote by $M \vee N$ the von Neumann algebra generated by $M$ and 
$N$.

\begin{lemma} Suppose $F_1$ and $F_2$ are finite sets of selfadjoint 
elements in $M$ and $x^* =x \in M$ satisfies $\chi(x) > -\infty$.  If 
$\{x\} 
\cup F_1$ and $\{x\} \cup F_2$ are $\alpha_1$-bounded and 
$\alpha_2$-bounded, then 
$\{x\} \cup F_1 \cup F_2$ is $(\alpha_1 + \alpha_2-1)$-bounded. 
\end{lemma}

\begin{proof} Fix $R$ greater than the operator norms of any of the 
elements in $F_1 \cup F_2$.  By Section 2 and Lemma 2.1 there exists a 
sequence $\langle x_k \rangle_{k=1}^{\infty}$ such that there exist 
constants $C_1, C_2, \epsilon_1, \epsilon_2 >0$ such that for any $m \in 
\mathbb N$ and $\gamma >0$, $x_k \in \Gamma(x;m,k,\gamma)$ for $k$ 
sufficiently large and such that the microstate spaces 
$\Xi(F_i;m,k,\gamma)$, $i=1,2$ relative to the sequence $\langle x_k 
\rangle_{k=1}^{\infty}$ satisfy for all $\epsilon_i >t >0$

\[ (\alpha_i -1 )  |\log t| + C_i \geq \mathbb K_t(\Xi(F_i)).
\]

\noindent Notice that we have arranged the same sequence $x_k$ with 
respect to which we consider the conditioned microstates for $F_1$ and 
$F_2$ (this is exactly what was arranged in Lemma 2.1 and Lemma 2.2).  
Because $\Xi(F_1;m,k,\gamma) \times 
\Xi(F_2;m,k,\gamma) \supset \Xi(F_1 \cup F_2;m,k,\gamma)$, it follows that 
for all $\min \{\epsilon_1, \epsilon_2\} > \epsilon >0$,

\begin{eqnarray*} (\alpha_1 + \alpha_2 -2) |\log \epsilon| + C_1 +C_2 & 
\geq & \mathbb K_{\epsilon}(\Xi(F_1)) + \mathbb K_{\epsilon} (\Xi(F_2)) \\ 
& \geq & \mathbb K_{\epsilon \sqrt{2}}(\Xi(F_1 \cup F_2)). \end{eqnarray*}

\noindent By Lemma 2.2 there exist $C, \epsilon_3 >0$ such that for all 
$\epsilon_3 > \epsilon >0$

\begin{eqnarray*} \mathbb K_{\epsilon, R}(\{x\} \cup F_1 \cup F_2) & \leq 
& C + |\log \epsilon| + \mathbb K_{\frac{\epsilon}{6}, R}(\Xi(F_1 \cup 
F_2)) \\ & \leq & C + C_1 + C_2 + 4(\alpha_1 +\alpha_2) + (\alpha_1 + 
\alpha_2 -1) |\log \epsilon|. \end{eqnarray*}

\noindent By Corollary 1.4 $\{x\} \cup F_1 \cup F_2$ is $(\alpha_1 + 
\alpha_2\ -1)$-bounded. \end{proof}

\begin{corollary} If $M$ and $N$ are strongly $1$-bounded and $M \cap N$ 
is diffuse, then $M \vee N$ is strongly one-bounded.  \end{corollary}

\begin{proof} By hypothesis, $M$ and $N$ can be generated by finite sets 
of selfadjoint elements $F$ and $G$, respectively, such that $F$ and $G$ 
each contain a selfadjoint with finite free entropy.  Pick a semicircular 
element $z \in M \cap N$.  By Theorem 3.2 $\{z\} \cup F$ is $1$-bounded as 
is $\{z\} \cup G$.  Thus, by Lemma 4.1, $\{z\} \cup F \cup G$ is a 
$1$-bounded set.  Clearly $\{z\} \cup F \cup G$ generates $M \vee N$ and 
$\chi(z) > -\infty$ by [14], so by definition $M \vee N$ is strongly 
one-bounded. \end{proof}

\begin{corollary} Suppose $M$ and $N$ are strongly $1$-bounded von Neumann 
algebras with diffuse von Neumann subalgebras $A$ and $B$, respectively.  
If $A \vee B$ is strongly $1$-bounded, then $M \vee N$ is strongly 
one-bounded.  In particular, $M \vee N$ is not isomorphic to an 
interpolated free group factor or to the free perturbation algebras in 
[1].  \end{corollary}

\begin{proof} By Corollary 4.2 the von Neumann algebras $M \vee B = M \vee 
(A \vee B)$ and $N \vee A = N \vee (B \vee A)$ are strongly $1$-bounded.  
Since $A \vee B$ is a diffuse von Neumann algebra contained in both $M 
\vee B$ and $N \vee A$, Corollary 4.2 implies that $M \vee N$ is strongly 
$1$-bounded.  The rest is given by Corollary 3.8. \end{proof}

\begin{corollary} If $M$ and $N$ are strongly $1$-bounded and $D \subset 
M\cap N$ is a diffuse subalgebra, then the amalgamated free product $M *_D 
N$ is not isomorphic to the interpolated free group factors or  
free perturbations algebras in [1]. \end{corollary}

\section{Normalizers and strongly $1$-bounded von Neumann algebras}

We now turn to a generalization of the von Neumann algebras in [5] and 
[17].  In this last section $A \subset M$ is an inclusion of tracial von 
Neumann algebras and $\{x\} \cup F$ is finite set of selfadjoint 
generators for $A$.  Assume $R >1$ exceeds the norms of any of the 
elements of $\{x\} \cup F$.  The relative microstates $\Xi()$ and 
associated quantities written below will all be taken with respect to a 
fixed sequence $\langle x_k \rangle_{k=1}^{\infty}$ of microstates for $x$ 
as discussed in Section 2.

\begin{lemma} Suppose that $u \in M$ is a unitary such that 
for some diffuse selfadjoint $y \in A$, $uyu^* \in A$, and 
$z \in M$ is a selfadjoint such that $z^{\prime \prime} = 
u^{\prime \prime}$.  If $1 > \epsilon, r >0$, then there 
exist $m_0 \in \mathbb N$, $\gamma_0 >0$ and a constant 
$L(\epsilon) >1$ dependent on $\epsilon$ such that for any 
$m > m_0$, $\gamma_0 > \gamma >0$, and $\epsilon 
L(\epsilon)^{-1}$-net $\langle \xi_{sk} \rangle_{s \in 
\Sigma_k}$ for $\Xi_R(F;m,k,\gamma)$, there exist an index 
set $\Theta_k$ satisfying $\#\Theta_k < \epsilon^{-rk^2}$ 
and for each $s \in S_k$ a collection $\langle \eta_{bsk} 
\rangle_{b \in \Theta_k}$ such that $\langle (\xi_{sk}, 
\eta_{bsk}) \rangle_{(s,b) \in \Sigma_k \times \Theta_k}$ is 
an $\epsilon$-cover for $\Xi_R(F \cup \{z\} ;m,k,\gamma)$.  
\end{lemma}

\begin{proof} Suppose $1 > \epsilon, r >0$.  There exists a polynomial $h$ 
in one $*$-variable such that $|h(u) -z |_2 < \epsilon (40)^{-1}$.  There 
also exists a constant $K >1$ so that regarding $h$ as a function from 
$(M_k(\mathbb C))_2$ into $M_k(\mathbb C)$, $h$ has a Lipschitz constant 
no greater than $K$, $K$ is independent of $k$.  
Find $n \in \mathbb N$ satisfying $nr |\log \epsilon| > \log (40K) + |\log 
\epsilon|$.  Choose mutually orthogonal projections $e_{1}, \ldots, e_{n} 
\in \{y\}^{\prime \prime}$, all with trace $n^{-1}$.

\[ u = \Sigma_{i=1}^n u e_i = \Sigma_{i=1}^n (u e_i u^*) u 
e_i.
\] 

\noindent For each $1 \leq i \leq n$ set $f_i = u e_i u_i^*$ 
and observe that $e_i, f_i \in A$.  Thus, $\Sigma_{i=1}^n 
f_i u e_i = u$ and $|h[\Sigma_{i=1}^n f_i u e_i] -z|_2 < 
\epsilon(10)^{-1}.$

From the first paragraph it follows that there exist 
noncommutative, selfadjoint polynomials $\Phi_i$, $\Psi_i$ 
in $\#F$-variables, $1 \leq i \leq n$, a polynomial $g$ in 
one variable, and $m_0 \in \mathbb N$, $\gamma_0 >0$ such 
that the following conditions are satisfied for any $m > 
m_0$, $\gamma_0 > \gamma >0$:

\begin{itemize}

\item If $(\xi, \eta) \in \Gamma_R(F \cup \{z\};m, k, 
\gamma)$, then $\left | h \left[ \Sigma_{i =1}^n \Psi_i(\xi) 
g(\xi) \Phi_i(\xi) \right] - \eta \right|_2 < 
\frac{\epsilon}{6}.$

\item For any $\#F$-tuple $\xi$ of selfadjoint elements in a von Neumann
algebra with entries with operator norm no greater than $R$ or any single
selfadjoint element $\eta$ in the von Neumann algebra with operator norm
no greater than $R$, $\Phi_i(\xi)$, $\Psi_i(\xi)$, and $g(\eta)$ are all
contractions.

\item If $\xi \in \Gamma(F;m, k, \gamma)$, then there exists 
projections $p_i, q_i \in M_k(\mathbb C)$ all with 
normalized trace $n^{-1}$ such that the $p_i$ are mutually orthogonal, the $q_i$ are mutually orthogonal, $|\Phi_i(\xi) - 
q_{ij}|_2 < \epsilon(40nK)^{-1}$, $|\Psi_i(\xi) - p_i|_2 < 
\epsilon(40nK)^{-1}$, and 

\[ \left | h \left[ \Sigma_{i =1}^n 
\Psi_i(\xi) g(\xi) \Phi_i(\xi) \right] - h \left[ \Sigma_{i 
=1}^n p_i g(\xi) q_i \right] \right|_2 < 
\epsilon(10K)^{-1}.\]

\end{itemize}

\noindent Observe that each of the $\Phi_i$ and $\Psi_i$, 
are Lipschitz when considered as maps from 
$(M^{sa}_k(\mathbb C)_R)^{\#F}$ into $M_k(\mathbb C)$ where 
the Lipschitz constants are independent of $k$ and the 
domains and ranges of these polynomial maps are endowed with 
the $| \cdot |_2$-norm. Thus, there exists a constant $C > 
1$ which exceeds the Lipschitz constants of any of the 
$\Phi_i$ or $\Psi_i$ so regarded.  Set $L(\epsilon) = 
40CKn$.

Now suppose for each $k$ $\langle \xi_{sk} \rangle_{s \in 
\Sigma_k}$ is an $\epsilon L(\epsilon)^{-1}$-net for 
$\Xi_R(F;m,k,\gamma)$.  For each $s \in \Sigma_k$ fix a 
family of projections $\langle p_{isk} \rangle_{i=1}^{n} 
\cup \langle q_{isk} \rangle_{i=1}^n$ all with normalized 
trace $n^{-1}$ such that the $p_{isk}$ are mutually orthogonal, the $q_{isk}$ are mutually orthogonal, $|\Phi_i(\xi_{sk}) - q_{isk}|_2 
< \epsilon (40nK)^{-1}$ and $|\Psi_i(\xi_{sk}) - 
p_{isk}|_2 < \epsilon (40nK)^{-1}$ (this is possible by the third condition of the second paragraph).  Consider  
$\Sigma_{i =1}^n q_{isk} (M^{sa}_k(\mathbb C))_1 
p_{isk} \subset (M^{sa}_k(\mathbb C))_1$. Find an 
$\epsilon(40K)^{-1}$-net $\langle \eta_{bs} \rangle_{b \in 
\Theta_{k}}$ for $\Sigma_{i=1}^n q_{isk} 
(M^{sa}_k(\mathbb C))_1 p_{isk}$ such that

\[ \# \Theta_{k} \leq \left(\frac{40K}{\epsilon} 
\right)^{\frac{k^2}{n}} \leq \left( 
\frac{1}{\epsilon} \right)^{rk^2} \]

\noindent Consider $\langle (\xi_{sk}, h(\eta_{bsk}) 
\rangle_{(s, b) \in \Sigma_k \times \Theta_k}$.  I claim 
that this is an $\epsilon$-net for $\Xi_R(F\cup \{z\} ; 
m,k,\gamma)$.  Towards this end suppose $(\xi, \eta) \in 
\Gamma_R(F \cup \{z\}; m,k,\gamma)$.  Denote by $p_i$ and 
$q_i$ the projections provided for in the third condition of 
the second paragraph.  There exists an $s_0 \in \Sigma_k$ 
such that $|\xi - \xi_{s_0k}|_2 < \epsilon 
L(\epsilon)^{-1}$. For any $1 \leq i \leq n$

\begin{eqnarray*} |p_i - p_{is_0k}|_2 & \leq & |p_i - 
\Psi_i(\xi)|_2 + |\Psi_i(\xi) - \Psi_i(\xi_{s_0k})|_2 + 
|\Psi_i(\xi_{s_0k}) - p_{is_0k}|_2 \\ & \leq & \epsilon 
(20nK)^{-1} + \epsilon (40nK)^{-1} = \epsilon 
(10nK)^{-1}.\\ \end{eqnarray*}

\noindent Similarly, $|q_i - q_{is_0k}|_2 < \epsilon 
(10nK)^{-1}$.  From this it follows that $|\Sigma_{i=1}^n 
p_i 
g(\eta) q_i - \Sigma_{i=1}^n p_{is_0k} g(\eta) q_{is_0k}|_2 
< \epsilon (5K)^{-1}$ whence, $|h(\Sigma_{i=1}^n p_i g(\eta) 
q_i) - h(\Sigma_{i=1}^n p_{is_0k} g(\eta) q_{is_0k})|_2 < 
\epsilon 5^{-1}$.  Since $\|g(\eta)\| \leq 1$, it follows 
that 
there exists some $b_0 \in \Theta_k$ satisfying

\[ \left |\Sigma_{i = 1}^n
q_{is_0k}f_i(\eta)p_{is_0k} - \eta_{b_0sk} \right |_2 <
\epsilon(40K)^{-1}.\]

\noindent Putting this all together we get

\begin{eqnarray*} |\eta - h(\eta_{bs_0k}) |_2 & \leq & |\eta 
-h [\Sigma_{i=1}^n \Psi_i(\xi) g(\eta) \Phi_i(\xi)]|_2 + \\ 
&& |h[\Sigma_{i=1}^n \Psi_i(\xi)g(\eta) \Phi_i(\xi)] - 
h[\Sigma_{i=1}^n p_i g(\eta) q_i]|_2 + \\ & & 
|h[\Sigma_{i=1}^n p_i g(\eta) q_i) - h[\Sigma_{i=1}^n 
p_{is_0k} g(\eta) q_{is_0 k}]|_2 + \\ && |h[\Sigma_{i=1}^n 
p_{is_0k} g(\eta) q_{is_0k}] - h(\eta_{b_0sk})|_2 \\ & \leq 
& \frac{\epsilon}{6} + \frac{\epsilon}{10} + 
\frac{\epsilon}{5} + \frac{\epsilon}{40} \\ & < & 
\frac{3\epsilon}{5}. 
\end{eqnarray*}

\noindent Thus, $|(\xi, \eta) - (\xi_{s_0k}, 
\eta_{b_0s_0k})|_2 < \epsilon L(\epsilon)^{-1} + 
\frac{3\epsilon}{5} < \epsilon$.  This completes the proof.
\end{proof}

Now suppose $\langle u_i \rangle_{i=1}^{\infty}$ is a 
sequence of unitaries in $M$.  For each $i$ let $z_i$ be a 
selfadjoint contraction such that $\{z_i\}^{\prime \prime} = 
\{u_i\}^{\prime \prime}$.  By scaling the $z_i$, we can 
arrange it so that $\sum_{i=1}^{\infty} \|z_i\|^2_2 <1$.  
Denote by $A_i$ the von Neumann algebra generated by $A$ and 
$\{u_1,\ldots, u_i\}$.  Assume that $u_1 \in A$ and for each 
each $i$, $u_{i+1}u_i u_{i+1}^* \in A_i$.

\begin{lemma} If $A$ is strongly $1$-bounded, then there 
exist $K,\epsilon_0 > 0$ such that for any $i \in \mathbb 
N$ and $\epsilon_0 > \epsilon >0$

\[ \mathbb K_{\epsilon, R}(\Xi(F \cup \{z_1,\ldots, z_i\})) < K.
\]
\end{lemma}

\begin{proof} We produce $C, \epsilon_0 >0$ such that for any $i \in N$ 
and $\epsilon_0 > \epsilon >0$,

\[ \mathbb K_{\epsilon}(\Xi_R(F \cup \{z_1, \ldots, z_i\})) \leq C + 
\Sigma_{j=1}^i 2^{-i}. \]

\noindent We will demonstrate this by induction on $i$.

Because $\{x\} \cup F$ is a $1$-bounded set and $\chi(x) > -\infty$ there 
exists by Lemma 2.1 constants $C, \epsilon_0 >0$, $1 > \epsilon_0$, 
dependent 
on $\{x\} \cup F$ 
such that for any $\epsilon_0 > \epsilon >0$,
 
\begin{eqnarray} C \geq \mathbb K_{\epsilon}(\Xi(F)). \end{eqnarray}

We can now start the induction.  Suppose $1 >\epsilon_0 > \epsilon >0$.  
Now $1 > \epsilon, 4^{-1} \epsilon >0$ so applying Lemma 5.1 with $r 
=4^{-1} \epsilon$ yields the corresponding $m \in \mathbb N$, $\gamma >0$, 
and constant $L(\epsilon) >1$ dependent on $\epsilon$ such that the 
conclusion of Lemma 5.1 holds.  $\epsilon L(\epsilon)^{-1} < \epsilon_0$ 
so (7) provides $m_1 >m$ and $\gamma > \gamma_1 >0$ so that $\mathbb 
K_{\epsilon L(\epsilon)^{-1}}(\Xi_R(F;m_1, \gamma_1)) \leq C + 4^{-1}$.  
Thus, for $k$ sufficiently large there exists an $\epsilon 
L(\epsilon)^{-1}$-net $\langle \xi_{sk} \rangle_{s \in \Sigma_k}$ for 
$\Xi_R(F;m_k,k,\gamma_1)$ such that $\# \Sigma_k < e^{k^2(C +4^{-1})}$.  
Applying the conclusion of Lemma 5.1 yields

\begin{eqnarray*} \mathbb K_{\epsilon}(\Xi_R(F \cup \{z_1\})) & \leq & 
\mathbb K_{\epsilon}(\Xi_R(F \cup \{z_1\};m_1, \gamma_1) \\ & \leq & C + 
4^{-1} + 4^{-1} \epsilon |\log \epsilon| \\ & \leq & C + 2^{-1}. 
\end{eqnarray*}

Now suppose the statement is true at $i \in \mathbb N$.  Suppose 
$\epsilon_0 > \epsilon >0$.  Again, applying Lemma 5.1 where $F$ is 
replaced by $F \cup \{z_1, \ldots, z_i\}$ $z$ is replaced by $z_{i+1}$, 
and $4 > \epsilon, 4^{-i}\epsilon =r >0$, there exists an $m \in \mathbb 
N$, $\gamma >0$, and $L(\epsilon) >1$ dependent on $\epsilon$ such that 
the conclusion of Lemma 5.1 holds.  By the inductive hypothesis,

\[ \mathbb K_{\epsilon L(\epsilon)^{-1}}(\Xi_R(F \cup \{z_1, \ldots, 
z_i\})) \leq C + \Sigma_{j=1}^i 2^{-j}. \]

\noindent Consequently, there exists $m_1 \in \mathbb N$ and $\gamma_1$ 
such that for $k$ sufficiently large there exists an $\epsilon 
L(\epsilon)^{-1}$-net $\langle \zeta_{nk} \rangle_{n \in \Omega_k}$ for 
$\Xi_R(F \cup \{x_1, \ldots, x_i\};m_1, k, \gamma_1)$ such that $\Omega_k 
< e^{k^2(C + 4^{-(i+1)} + \Sigma_{j=1}^i 2^{-j})}$.  Now applying the 
conclusion of Lemma 5.1 to these nets, it follows that

\begin{eqnarray*} \mathbb K_{\epsilon}(\Xi_R(F \cup \{z_1, \ldots, 
z_{i+1}\})) & \leq & \mathbb K_{\epsilon}(\Xi_R(F \cup \{z_1, \ldots, 
z_{i+1}\};m_1, \gamma_1) \\ & \leq & C + \Sigma_{j=1}^i 2^{-j} + 
4^{-(i+1)} + 4^{-(i+1)} \epsilon |\log \epsilon| \\ & < & C + 
\sum_{j=1}^{i+1} 2^{-j}. \end{eqnarray*} \end{proof}

We are now ready for the main result of the section.  Its proof runs very 
much like that of Corollary 3.2, except that issues with normalizers 
complicate the argument (hence the preceding prepatory lemmas).  For 
efficiency's sake, we could have subsumed Theorem 3.2 with the lemmas and 
the following general theorem of this section.  But the relation between 
the invariance issue and strong $1$-boundedness is clearer in the less 
cluttered context of Theorem 3.2, and for clarity's sake (clarity and 
efficiency not being the same) we have opted to repeat (more or less) the 
argument of Theorem 3.2 below.

\begin{theorem} If $A_{\infty}$ is the von Neumann algebra generated by 
$A$ and $\langle u_i \rangle_{i=1}^{\infty}$, then any finite set of 
selfadjoint generators $G$ for $A_{\infty}$ is $1$-bounded. \end{theorem}

\begin{proof} Let $K, \epsilon_0$ be as in the preceding lemma.  Without 
loss of generality we can assume the operator norms of any of the elements 
in $G$ is no greater than $R$.  Recall the universal constant $L$ in Lemma 
2.2 and find a $D > 4[\Sigma_{y \in G} \|y\|_2^2 + \|x\|_2^2 + 
\Sigma_{i=1}^{\infty} \|z_i\|_2^2 +6] +L$.  By Lemma 3.1 and Lemma 2.2 it 
follows that for any $i \in \mathbb N$ and $ \epsilon_0 > \epsilon >0$,

\begin{eqnarray*} \mathbb K_{\epsilon, R}(G) & \leq &  \mathbb K_{\epsilon, R}( \{x\} \cup F \cup \{z_1, \ldots z_i\} \cup G) \\ & \leq & \log D + |\log \epsilon| + \mathbb K_{\epsilon D^{-1}, R}(\Xi(F \cup \{z_1,\ldots, z_i\} \cup G)) \\
\end{eqnarray*}

\noindent To finish the proof then, it suffices to bound the last term on 
the dominating sum above.  This will follow from arguing just as in 
Theorem 3.2.

Fix $\epsilon$ with $\epsilon_0 > \epsilon > 0$.  There exists an $N$ so 
large and a $\#G$-tuple of noncommuting polynomials $\Phi$ in $(\#F + N 
+1)$-variables satisfying $|\Phi(x, F, z_1, \ldots, z_n) - G|_2 < \epsilon 
(4D)^{-1}$.  Regarding $\Phi$ as a function from $((M^{sa}_k(\mathbb 
C)_R))^{\#F + N+1}$ into $(M^{sa}_k(\mathbb C))^{\#G}$, $\Phi$ has bounded 
Lipschitz constant $L_k$ with respect to the $|\cdot |_2$-norm and 
moreover, $\sup_{k \in \mathbb N} L_k < \infty$.  Choose $L > \sup_{k \in 
\mathbb N} L_k + 1$.  Suppose $t$ satisfies $\epsilon(4DL)^{-1} > t >0$.  
By Lemma 5.2

\[ \mathbb K_{t, R}(\Xi(F \cup \{z_1, \ldots, z_n\})) <  K.\]

Thus, there exist $m_1 \in \mathbb N$ and $\gamma_1 >0$ such that for $m 
\geq m_1$, and $\gamma_1 > \gamma >0$,

\begin{eqnarray} K+1 \geq \limsup_{k \rightarrow \infty} k^{-2} \cdot \log 
K_t(\Xi_R(F \cup \{z_1, \ldots, z_n\};m,k,\gamma)). \end{eqnarray}

There exist $m_2 \in \mathbb N$ and $\gamma_2 >0$ such that if $m > m_2$ 
and $\gamma_2 > \gamma >0$, then for any $k, (\xi, \eta, \zeta) \in 
\Xi(\{x\} \cup F \cup \{z_1,\ldots, z_n\} \cup G;m, k, \gamma) \Rightarrow 
|\Phi(x_k, \xi, \eta) - \zeta|_2 \leq \epsilon (4D)^{-1}$. Suppose $m > 
m_1 + m_2$ and $\min\{\gamma_1, \gamma_2\} > \gamma >0$.  By (8) for $k$ 
sufficiently large we can find a $t$-net $\langle (\xi_{sk}, \eta_{sk}) 
\rangle_{s \in \Sigma_k}$ for $\Xi_R(F \cup \{x_1, \ldots, 
x_n\};m,k,\gamma)$ which satisfies

\begin{eqnarray} \# \Sigma_k \leq e^{(K+1)k^2}.
\end{eqnarray}

\noindent For each such $k$ sufficiently large consider $\langle 
(\xi_{sk}, \eta_{sk}, \Phi(x_k, \xi_{sk}, \eta_{sk})) \rangle_{s \in 
\Sigma_k}$.  This set is an $\epsilon D^{-1}$-cover for $\Xi_R(F \cup 
\{z_1, \ldots, z_n\} \cup G;m,k,\gamma)$.  To see this, suppose $(\xi, 
\eta, \zeta) \in \Xi_R(F \cup \{z_1,\ldots, z_n\} \cup G;m, k,\gamma)$.  
By definition $(\xi, \eta) \in \Xi_R( F \cup \{z_1,\ldots, 
z_n\};m,k,\gamma)$ so there exists some $s_0 \in \Sigma_k$ with 
$|(\xi_{s_0k}, \eta_{s_0k}) - (\xi, \eta)|_2 \leq t \leq 
\epsilon(4DL)^{-1}$.  Both $(\eta, \zeta)$ and $(\eta_{s_0k}, 
\zeta_{s_0k})$ 
are $(\#F+n)$-tuples of operators with norms no greater than $R$ so 
$|\Phi(x_k, \xi, \eta) - \Phi(x_k, \xi_{s_0k}, \eta_{s_0k})|_2 \leq L 
\cdot t \leq \epsilon (4D)^{-1}$.  But $|\Phi(x_k, \xi, \eta) -\zeta|_2 < 
\epsilon (4D)^{-1}$ so $|\Phi(x_k, \xi_{s_0k}, \eta_{s_0k}) - \zeta|_2 < 
\epsilon(2D)^{-1}$.  Combining this with $|(\xi_{s_0k}, \eta_{s_0k}) 
-(\xi, \eta)|_2 <t$ yields

\[ |(\xi, \eta, \zeta) - (\xi_{s_0k}, \eta_{s_0k}, \Phi(x_k, 
\xi_{s_0k}, \eta_{s_0k}))|_2 < 
\epsilon D^{-1}.\]

It now follows that for $m > m_1 + m_2$ and $\min\{\gamma_1,\gamma_2\} > 
\gamma >0$ and $k$ sufficiently large.

\[ k^{-2} \cdot \log \left [ K_{\epsilon D^{-1}}(\Xi_R(F \cup 
\{z_1,\ldots,z_n\} \cup G;m,k,\gamma)) \right] \leq K+1. \]

\noindent This implies

\begin{eqnarray} \mathbb K_{\epsilon D^{-1}, R}(\Xi(F \cup \{z_1,\ldots, 
z_n\} \cup G)) \leq K+1. \end{eqnarray}

\noindent $K$ is dependent only on $\{x\} \cup F$ and thus, by the first 
paragraph, we're done. \end{proof}

In Theorem 5.3 if $A_{\infty}$ is a $\mathrm{II}_1$-factor, 
then a trivial modification of the techniques in [6] shows that 
$A_{\infty}$ can be generated by a finite set of selfadjoint elements.  By 
Theorem 
5.3 then, $A_{\infty}$ is strongly $1$-bounded.  We thus have,

\begin{corollary} Suppose $A$ is a strongly $1$-bounded von Neumann 
algebra and $\langle u_i \rangle_{i=1}^{\infty}$ is a sequence of Haar 
unitaries such that $u_1 \in A$ and for each $i \in \mathbb N$ $u_{i+1} 
u_i u_{i+1}^* 
\in (A \cup \{u_1, \ldots u_i \})^{\prime \prime}$.  If the von 
Neumann algebra $A_{\infty}$ generated by $A$ and $\langle u_i 
\rangle_{i=1}^{\infty}$ is a $\mathrm{II}_1$-factor, then $A_{\infty}$ is 
not isomorphic to an interpolated free group factor, $L(F_r)$, $1 < r \leq 
\infty$.  In particular, a $\mathrm{II}_1$-factor generated by the 
normalizer of a strongly $1$-bounded von Neumann subalgebra is not an 
interpolated free group factor. \end{corollary}

\begin{remark} It is immediate from [7] that any diffuse, hyperfinite von 
Neumann algebra is strongly $1$-bounded.  Thus, Theorem 5.3 provides an 
alternate way of seeing that the von Neumann algebras considered in [5] or 
[17] have free entropy dimension no greater than $1$.  \end{remark}

\end{document}